\overfullrule=0pt
\centerline {\bf An improvement of a saddle point theorem and some of its applications}\par
\bigskip
\bigskip
\centerline {BIAGIO RICCERI}\par
\bigskip
\bigskip
\centerline {\it Dedicated to the memory of Professor Wataru Takahashi}\par
\bigskip
\bigskip
\bigskip
\bigskip
{\bf Abstract.} In this paper, we establish an improved version of a saddle point theorem ([4]) removing a weak lower semicontinuity assumption at all. We then revisit some of the applications of that theorem in the light of such an improvement. For instance, we obtain the following very general result of local nature: Let $(H,\langle\cdot,\cdot\rangle)$ be a real Hilbert space
and $\Phi:B_{\rho}\to H$ a $C^{1,1}$ function, with $\Phi(0)\neq 0$. Then, for each $r>0$ small enough, there exist only two points points $x^*, u^*\in S_r$,
such that
$$\max\{\langle \Phi(x^*),x^*-x\rangle, \langle \Phi(x),x^*-x\rangle\}< 0$$
for all $x\in B_r\setminus \{x^*\}$,
$$\|\Phi(u^*)-u^*\|=\hbox {\rm dist}(\Phi(u^*),B_r)$$
and
$$\|\Phi(x)-u^*\|<\|\Phi(x)-x\|$$
for all $x\in B_r\setminus \{u^*\}$, where
$$B_r=\{x\in H : \|x\|\leq r\}$$
and
$$S_r=\{x\in H : \|x\|=r\}\ .$$
\par
\bigskip
\bigskip
\bigskip
\bigskip
In the sequel, $(H,\langle\cdot,\cdot\rangle)$ is a real Hilbert space. For each $r>0$, set
$$B_r=\{x\in H : \|x\|\leq r\}$$
and
$$S_r=\{x\in H : \|x\|=r\}\ .$$
The aim of this paper is to establish the following result jointly with two meaningful applications of it:\par
\medskip
THEOREM 1. - {\it Let $Y$ be a non-empty closed convex set in a Hausdorff real topological vector space, let $\rho>0$ and
let $J:B_{\rho}\times Y\to {\bf R}$ be a function satisfying the following conditions:\par
\noindent
$(a_1)$\hskip 5pt  for each $y\in Y$, the function $J(\cdot,y)$ is $C^1$ and $J'_x(\cdot,y)$ is Lipschitzian
with constant $L$ (independent of $y$)\ ;\par
\noindent
$(a_2)$\hskip 5pt $J(x,\cdot)$ is upper semicontinuous and concave for all $x\in B_{\rho}$ and $J(x_0,\cdot)$ is
sup-compact for some $x_0\in B_{\rho}$;\par
\noindent
$(a_3)$\hskip 5pt $\delta:=\inf_{y\in Y}\|J'_x(0,y)\|>0\ .$\par
Then, for each $r\in \left ] 0,\min\left \{\rho, {{\delta}\over {2L}}\right \}\right ]$ and for each non-empty closed convex $T\subseteq Y$, there exist
$x^*\in S_r$ and $y^*\in T$ such that
$$J(x^*,y)\leq J(x^*,y^*)< J(x,y^*)$$
for all $x\in B_r\setminus \{x^*\}$, $y\in T$\ .}\par
\smallskip
PROOF. Fix $r\in \left ] 0,\min\left \{\rho, {{\delta}\over {2L}}\right \}\right ]$ and a non-empty closed convex $T\subseteq Y$. Consider the function
$\varphi:B_r\times T\to {\bf R}$ defined by
$$\varphi(x,y)={{L}\over {2}}\|x\|^2+J(x,y)$$
for all $(x,y)\in B_r\times T$. Notice that, for each $y\in T$, the function $\varphi(\cdot,y)$ is continuous and convex in $B_r$ (see the proof of Corollary 2.3 of [3]).
Consequently, if we consider $B_r$ endowed with the relative weak topology, the function $\varphi$ satisfies the assumptions of a classical minimax theorem ([1], Theorem 2) from which
we infer
$$\sup_T\inf_{B_r}\varphi=\inf_{B_r}\sup_T\varphi\ .\eqno{(1)}$$
The function $x\to \sup_{y\in T}\varphi(x,y)$ (resp. $y\to \inf_{x\in B_r}\varphi(x,y)$) is weakly lower semicontinuous (resp. sup-compact).  Therefore, there exist $x^*\in B_r$ and $y^*\in T$ such that
$$\sup_{y\in T}\varphi(x^*,y)=\inf_{x\in B_r}\sup_{y\in T}\varphi(x,y)\ ,$$
$$\inf_{x\in B_r}\varphi(x,y^*)=\sup_{y\in T}\inf_{x\in B_r}\varphi(x,y)\ .$$
So, in view of $(1)$, we obtain
$$\varphi(x^*,y)\leq \varphi(x^*,y^*)\leq \varphi(x,y^*)\eqno{(2)}$$
for all $x\in B_r$, $y\in T$. 
Notice that the equation
$$J'_x(x,y^*)+L x=0$$
has no solution in the interior of $B_r$. Indeed, let $\tilde x\in B_{\rho}$ be such that
$$J'_x(\tilde x,y^*)+L\tilde x=0\ .$$
Then, in view of $(a_1)$, we have
$$\|L\tilde x+J'_x(0,y^*)\|\leq \|L\tilde x\|\ .$$
In turn, using the Cauchy-Schwarz inequality, this readily implies that
$$\|\tilde x\|\geq {{\|J'_x(0,y^*)\|}\over {2L}}\geq {{\delta}\over {2L}}\geq r\ .$$
From this remark, we infer that the set of all global minima of the function $\varphi(\cdot,y^*)$ is contained in $S_r$ and so, being
convex, it reduces to $x^*$ (recall that $X$ is a Hilbert space), in view of $(2)$.
 Therefore, for every $x\in B_r\setminus \{x^*\}$, $y\in T$, from $(2)$ we obtain
$${{1}\over {2}}\|x^*\|^2+J(x^*,y)\leq
{{1}\over {2}}\|x^*\|^2+J(x^*,y^*)<{{1}\over {2}}\|x\|^2+J(x,y^*)\leq {{1}\over {2}}\|x^*\|^2+J(x,y^*)$$
and so
$$J(x^*,y)\leq J(x^*,y^*)<J(x,y^*)\ .$$
The proof is complete.\hfill $\bigtriangleup$\par
\medskip
REMARK 1. - Theorem 1 was obtained in [4] ([4], Theorem 2.1) under the following additional assumption: for each $y\in Y$, the function $J(\cdot,y)$ is weakly lower semicontinuous. This was due to the
fact that, instead of applying the classical minimax theorem in [1] to $\varphi$ (as we did above), we applied Theorem 1.2 of [2] to $J$.\par
\medskip
We can now revisit two applications of Theorem 1. The first one concerns variational inequalities.\par
\medskip
THEOREM 2. - {\it Let $\rho>0$ and let $\Phi:B_{\rho}\to H$ be a $C^1$ function whose derivative is Lipschitzian with
constant $\gamma$. Set
$$\theta:=\sup_{x\in B_{\rho}}\|\Phi'(x)\|_{{\cal L}(H)}\ ,$$
$$M:=2(\theta+\rho\gamma)$$
and assume also that
$$\sigma:=\inf_{y\in B_{\rho}}\sup_{\|u\|=1}|\langle\Phi(0),u\rangle - \langle\Phi'(0)(u),y\rangle|>0\ .$$
Then, for each $r\in \left ] 0,\min\left \{\rho, {{\sigma}\over {2M}}\right \}\right ]$,  there exists a unique $x^*\in S_r$
such that 
$$\max\{\langle \Phi(x^*),x^*-x\rangle, \langle \Phi(x),x^*-x\rangle\}< 0$$
for all $x\in B_r\setminus \{x^*\}$.}\par
\smallskip
PROOF. Consider the function $J:B_{\rho}\times B_{\rho}\to {\bf R}$ defined by
$$J(x,y)=\langle\Phi(x),x-y\rangle$$
for all $x, y\in B_{\rho}$. Of course, for each $y\in B_{\rho}$, the function $J(\cdot,y)$ is $C^1$ and one has
$$\langle J'_x(x,y),u\rangle = \langle \Phi'(x)(u),x-y\rangle + \langle\Phi(x),u\rangle$$
for all $x\in B_{\rho}, u\in H$. Fix $x, v\in B_{\rho}$ and $u\in S_1$. We then have
$$|\langle J'_x(x,y),u\rangle - \langle J'_x(v,y),u\rangle|=|\langle\Phi(x)-\Phi(v),u\rangle +\langle \Phi'(x)(u),x-y\rangle-
\langle\Phi'(v)(u),v-y\rangle|$$
$$\leq \|\Phi(x)-\Phi(v)\|+|\Phi'(x)(u)-\Phi'(v)(u),v-y\rangle + \langle\Phi'(x)(u),x-v\rangle|$$
$$\leq\theta\|x-v\|+2\rho\|\Phi'(x)-\Phi'(v)\|_{{\cal L}(H)}+\theta\|x-v\|$$
$$\leq 2(\theta + \rho\gamma)\|x-v\|\ .$$
Hence, the function $J'_x(\cdot,y)$ is Lipschitzian with constant $M$. At this point, we can apply Theorem 1 taking $Y=B_{\rho}$ with
the weak topology. Therefore, for each $r\in \left ] 0,\min\left \{\rho, {{\sigma}\over {2M}}\right \}\right ]$,
there exist $x^*\in S_r$ and $ y^*\in B_{r}$ such that
$$\langle \Phi(x^*),x^*-y\rangle\leq \langle\Phi(x^*),x^*-y^*\rangle<\langle\Phi(x),x-y^*\rangle \eqno{(3)}$$
for all $x, y\in B_r$, with $x\neq x^*$. 
Notice that $\Phi(x^*)\neq 0$. Indeed, if $\Phi(x^*)=0$, we would have
$$\|\Phi(0)\|=\|\Phi(0)-\Phi(x^*)\|\leq \theta r$$
and hence, since $\sigma\leq \|\Phi(0)\|$, it would follow that
$$r\leq {{\|\Phi(0)\|}\over {2M}}<{{\|\Phi(0)\|}\over {\theta}}\leq r\ .$$ 
From the first inequality in $(3)$, taking $y=x^*$, we get $0\leq \langle\Phi(x^*),x^*-y^*\rangle$. So, in view of the
strict inequality, we infer that $x^*=y^*$ (since, otherwise, we could take $x=y^*$, obtaing a contradiction). Thus, $(3)$ actually reads
$$\langle \Phi(x^*),x^*-y\rangle\leq 0<\langle\Phi(x),x-x^*\rangle$$
for all $x, y\in B_r$, with $x\neq x^*$. In particular, we infer that $x^*$ is the unique global minimum in $B_r$ of the linear functional
$y\to \langle\Phi(x^*),y\rangle$.
Hence 
$$\max\{\langle \Phi(x^*),x^*-x\rangle, \langle \Phi(x),x^*-x\rangle\}< 0$$
for all $x\in B_r\setminus \{x^*\}$. Finally, to show the uniqueness of $x^*$, argue by contradiction, supposing that there is another
$\tilde x\in S_r$, with $\tilde x\neq x^*$, such that 
$$\max\{\langle \Phi(\tilde x),\tilde x-x\rangle, \langle \Phi(x)),\tilde x-x\rangle\}< 0$$
for all $x\in B_r\setminus \{\tilde x\}$. So, we would have at the same time 
$$\langle\Phi(\tilde x),\tilde x-x^*\rangle<0$$
and
$$\langle\Phi(\tilde x),x^*-\tilde x\rangle<0\ ,$$
an absurd, and the proof is complete.\hfill $\bigtriangleup$\par
\medskip
REMARK 2. - Theorem 2 was obtained in [4] ([4], Theorem 2.2) under the following additional assumption: for each $y\in B_{\rho}$, the
function $x\to \langle \Phi(x),x-y\rangle$ is weakly lower semicontinuous.\par
\medskip
A remarkable corollary of Theorem 2 is as follows:\par
\medskip
THEOREM 3. - {\it Let $\rho>0$ and let $\Phi:B_{\rho}\to H$ be a $C^1$ function with Lipschitzian derivative.\par
Then, the following assertions are equivalent:\par
\noindent
$(i)$\hskip 5pt for each $r>0$ small enough, there exists a unique $x^*\in S_r$ such that
$$\max\{\langle \Phi(x^*),x^*-x\rangle, \langle \Phi(x),x^*-x\rangle\}< 0$$
for all $x\in B_r\setminus \{x^*\}$\ ;\par
\noindent
$(ii)$\hskip 5pt $\Phi(0)\neq 0$\ .}\par
\smallskip
PROOF. The implication $(ii)\to (i)$ is obvious. So, assume that $(i)$ holds. Notice that the function
$y\to \sup_{\|u\|=1}|\langle\Phi(0),u\rangle - \langle\Phi'(0)(u),y\rangle|$ is continuous in $H$ and
takes the value $\|\Phi(0)\|>0$ at $0$. Therefore, for a suitable $r^*\in ]0,\rho]$, we have
$$\inf_{y\in B_{r^*}}\sup_{\|u\|=1}\langle\Phi(0),u\rangle - \langle\Phi'(0)(u),y\rangle|>0\ .$$
Now, we can apply Theorem 2 to the restriction of the function $\Phi$ to $B_{r^*}$, and $(i)$ follows.
\hfill $\bigtriangleup$\par
\medskip
Also, it is worth noticing the following further corollary of Theorem 2:\par
\medskip
THEOREM 4. - {\it Let $\rho>0$ and let $\Psi:B_{\rho}\to H$ be a $C^1$ function whose derivative vanishes at $0$ and is Lipschitzian with
constant $\gamma_1$. Set
$$\theta_1:=\sup_{x\in B_{\rho}}\|\Psi'(x)\|_{{\cal L}(H)}\ ,$$
$$M_1:=2(\theta_1+\rho\gamma_1)$$
and let $w\in H$ satisfy
$$\|w-\Psi(0)\|\geq 2M_1\rho\ .\eqno{(4)}$$
Then, for each $r\in ] 0,\rho]$,  there exists a unique $x^*\in S_r$
such that 
$$\max\{\langle \Psi(x^*)-w,x^*-y\rangle, \langle \Psi(y)-w,x^*-y\rangle\}< 0$$
for all $y\in B_r\setminus \{x^*\}$.}\par
\smallskip
PROOF. Set $\Phi:=\Psi-w$. Apply Theorem 2 to $\Phi$. Since $\Phi'=\Psi'$, we have $M=M_1$. Since $\Phi'(0)=0$, we have
$\sigma=\|\Phi(0)\|$ and so, in view of $(4)$, 
$$\rho\leq {{\sigma}\over {2M}}$$
and the conclusion follows.\hfill $\bigtriangleup$\par
\medskip
The second application of Theorem 1 is as follows:\par
\medskip
THEOREM 5. - {\it Let $Y\subseteq H$ be a closed bounded convex set, let $\rho>0$ and let $f:B_{\rho}\to H$ be a $C^1$ function whose derivative is Lipschitzian with constant $\gamma$. Moreover,
let $\eta$ be the Lipschitz constant of the function $x\to x-f(x)$, set
$$\theta:=\sup_{x\in B_{\rho}}\|f'(x)\|_{{\cal L}(H)}\ ,$$
$$L:=2\left (\eta+\theta+\gamma\left (\rho+\sup_{y\in Y}\|y\|\right )
\right )$$
and assume that
$$\sigma:=\inf_{y\in Y}\sup_{\|u\|=1}|\langle f'(0)(u),y\rangle-\langle f(0),u\rangle|>0\ .$$
Then, for each $r\in \left ] 0,\min\left \{\rho, {{\sigma}\over {L}}\right \}\right]$ and for each non-empty closed convex set $T\subseteq Y$, there exist
$x^*\in S_r$ and $y^*\in T$ such that
$$\|x^*-f(x^*)\|^2+\|f(x)-y^*\|^2-\|x-f(x)\|^2< \|f(x^*)-y^*\|^2=(\hbox {\rm dist}(f(x^*),T))^2 \eqno{(5)}$$
for all $x\in B_r\setminus \{x^*\}$\ .}\par
\smallskip
PROOF. Consider the function $J:B_{\rho}\times Y\to {\bf R}$ defined by
$$J(x,y)=\|f(x)-x\|^2-\|f(x)-y\|^2$$
for all $x\in B_{\rho}$, $y\in Y$. Clearly, for each $y\in Y$, $J(\cdot,y)$ is of class $C^1$. Moreover, one has
$$\langle J'_x(x,y),u\rangle = 2\langle x-f(x),u\rangle - 2\langle f'(x)(u),x-y\rangle$$
for all $x\in B_{\rho}$, $u\in H$. Fix $x, v\in B_{\rho}$ and $u\in H$, with $\|u\|=1$. We have
$${{1}\over {2}}|\langle J'_x(x,y)-J'_x(v,y),u\rangle|=|\langle x-f(x)-v+f(v),u\rangle-\langle f'(x)(u),x-y\rangle+\langle f'(v)(u),v-y\rangle|$$
$$\leq \eta\|x-v\|+|\langle f'(x)(u),x-v\rangle+\langle f'(x)(u)-f'(v)(u),v-y\rangle|$$
$$\leq \eta\|x-v\|+\|f'(x)(u)\|\|x-v\| + \|f'(x)(u)-f'(v)(u)\|\|v-y\|\leq \left (\eta+\theta+\gamma\left ( \rho+\sup_{y\in Y}\|y\|\right )\right )\|x-v\|\ .$$
Therefore, the function $J_x'(\cdot,y)$ is Lipschitzian with constant $L$. Moreover, we clearly have
$$\inf_{y\in Y}\|J'_x(0,y)\|=2\sigma\ .$$
So, if we fix $r\in \left ] 0,\min\left \{\rho, {{\sigma}\over {L}}\right \}\right ]$ and  a non-empty closed convex set $T\subseteq Y$, by Theorem 1,
there exist $x^*\in B_r$ and $y^*\in T$ such that
$$\|f(x^*)-x^*\|^2-\|f(x^*)-y\|^2\leq \|f(x^*)-x^*\|^2-\|f(x^*)-y^*\|^2<\|f(x)-x\|^2-\|f(x)-y^*\|^2 \eqno{(6)}$$
for all $x\in B_r\setminus \{x^*\}$, $y\in T$. Clearly, $(6)$ is equivalent to $(5)$, and the proof is complete.
\hfill $\bigtriangleup$
\medskip
REMARK 2. - Theorem 5 was obtained in [3] ([3], Corollary 2.5) assuming, in addition, that $f$ is sequentially weakly-strongly continuous.\par
\medskip
Here is a remarkable consequence of Theorem 5.\par
\medskip
THEOREM 6. - {\it Let $\rho>0$ and let $f:B_{\rho}\to H$ be a $C^1$ function whose derivative is Lipschitzian with constant $\gamma$. Moreover,
let $\eta$ be the Lipschitz constant of the function $x\to x-f(x)$, set
$$\theta:=\sup_{x\in B_{\rho}}\|f'(x)\|_{{\cal L}(H)}\ ,$$
$$L:=2(\eta+\theta+2\gamma\rho)$$
and assume that
$$\sigma:=\inf_{y\in B_{\rho}}\sup_{\|u\|=1}|\langle f'(0)(u),y\rangle-\langle f(0),u\rangle|>0\ .$$
Then, for each $r\in \left ] 0,\min\left \{\rho, {{\sigma}\over {L}}\right \}\right]$, there exists a unique $x^*\in S_r$ such that
$$\|f(x^*)-x^*\|=\hbox {\rm dist}(f(x^*),B_r) \eqno{(7)}$$
and
$$\|f(x)-x^*\|<\|f(x)-x\| \eqno{(8)}$$
for all $x\in B_r\setminus \{x^*\}$}.\par
\smallskip
PROOF. Fix $r\in \left ] 0,\min\left \{\rho, {{\sigma}\over {L}}\right \}\right]$. Applying Theorem 5  with $Y=B_{\rho}$ and $T=B_r$, we obtain $x^*\in S_r$ and $y^*\in B_r$
such that
$$\|x^*-f(x^*)\|^2+\|f(x)-y^*\|^2-\|x-f(x)\|^2< \|f(x^*)-y^*\|^2=(\hbox {\rm dist}(f(x^*),B_r))^2 \eqno{(9)}$$
for all $x\in B_r\setminus \{x^*\}$\ . From this, we infer that $y^*= x^*$. Actually, if $y^*\neq x^*$, we could take $x=y^*$ in
$(9)$, obtaining
$$\|x^*-f(x^*)\|<\hbox {\rm dist}(f(x^*),B_r)$$
which is absurd. Now, $(7)$ and $(8)$ follow directly from $(9)$. Finally, concerning the uniqueness of $x^*$, assume that $x_0\in B_r$
is such that
$$\|f(x)-x_0\|<\|f(x)-x\|$$
for all $x\in B_r\setminus \{x_0\}$. Then, if $x_0\neq x^*$, we would have
$$\|f(x^*)-x_0\|<\|f(x^*)-x^*\|$$
which is absurd in view of $(7)$.\hfill $\bigtriangleup$\par
\medskip
Finally, reasoning as in the proof of Theorem 3, we get the following corollary of Theorem 6:\par
\medskip
THEOREM 7. - {\it Let $\rho>0$ and let $f:B_{\rho}\to H$ be a $C^1$ function with Lipschitzian derivative.\par
Then, the following assertions are equivalent:\par
\noindent
$(i)$\hskip 5pt for each $r>0$ small enough, there exists a unique $x^*\in S_r$ such that
$$\|f(x^*)-x^*\|=\hbox {\rm dist}(f(x^*),B_r)$$
and
$$\|f(x)-x^*\|<\|f(x)-x\|$$
for all $x\in B_r\setminus \{x^*\}$\ ;\par
\noindent
$(ii)$\hskip 5pt $f(0)\neq 0$\ .}\par
\bigskip
\bigskip
{\bf Acknowledgement.} The author has been supported by the Universit\`a di Catania, PIACERI 2020-2022, Linea di intervento 2, Progetto “MAFANE” and by the Gruppo Nazionale per l'Analisi Matematica, la Probabilit\`a
e le loro Applicazioni (GNAMPA) of the Istituto Nazionale di Alta Matematica (INdAM).

\vfill\eject
\centerline {\bf References}\par
\bigskip
\bigskip
\noindent
[1]\hskip 5pt K. FAN, {\it Minimax theorems}, Proc. Nat. Acad. Sci. U.S.A., {\bf 39} (1953), 42-47.\par
\smallskip
\noindent
[2]\hskip 5pt B. RICCERI, {\it On a minimax theorem: an improvement, a new proof and an overview of its applications},
Minimax Theory Appl., {\bf 2} (2017), 99-152.\par
\smallskip
\noindent
[3]\hskip 5pt B. RICCERI, {\it Applying twice a minimax theorem}, J. Nonlinear Convex Anal., {\bf 20} (2019), 1987-1993.\par
\smallskip
\noindent
[4]\hskip 5pt B. RICCERI, {\it A remark on variational inequalities in small balls}, J. Nonlinear Var. Anal., {\bf 4} (2020), 21-26.\par
\bigskip
\bigskip
\bigskip
\bigskip
Department of Mathematics and Informatics\par
University of Catania\par
Viale A. Doria 6\par
95125 Catania, Italy\par
{\it e-mail address}: ricceri@dmi.unict.it

\bye

\bye